\documentclass[12pt]{article}
\usepackage{amsmath, amssymb}

\title{Meromorphic groups.}
\author{Anand Pillay\thanks{Supported by an NSF grant}\\Dept.
Mathematics\\ Univ. Illinois at Urbana-Champaign \and Thomas
Scanlon \thanks{Partially supported by an NSF MSPRF}
\\Dept. Mathematics\\University of California at
Berkeley}
\newtheorem{Theorem}{Theorem}[section]
\newtheorem{Proposition}[Theorem]{Proposition}
\newtheorem{Definition}[Theorem]{Definition}
\newtheorem{Lemma}[Theorem]{Lemma}
\newtheorem{Corollary}[Theorem]{Corollary}
\newtheorem{Remark}[Theorem]{Remark}
\newtheorem{Fact}[Theorem]{Fact}
\newtheorem{Conjecture}[Theorem]{Conjecture}

\newenvironment{proof}{\bigbreak \noindent {\bf Proof:}}{\hspace{\fill}
$\dashv$ \bigbreak}

\begin{document}
\maketitle
\begin{abstract} We introduce the notion of a meromorphic
group, weakening somewhat Fujiki's definition
(\cite{Fujiki}). We prove that a meromorphic group is
meromorphically an extension of a complex torus by a linear
algebraic group, generalizing results in \cite{Fujiki}. A
special case of this result, as well as one of the
ingredients in the proof, is that a strongly minimal
``modular" meromorphic group is a complex torus, answering a
question of Hrushovski. As a consequence, we show that
a simple compact complex manifold has algebraic and Kummer
dimension zero if and only if its generic type is trivial.
\end{abstract}

\section{Introduction}
Let ${\mathcal A}$ be the category of reduced irreducible
compact complex spaces. By a Zariski open subset of some
$X \in {\mathcal A}$ we mean, as usual, the complement of an
analytic subset of $X$. We understand a meromorphic mapping
from $X$ to
$Y$ ($X,Y\in {\mathcal A}$) in the sense of Remmert. Roughly
speaking it is an analytic subset $Z$ of $X\times Y$, such
that for some (nonempty, so dense) Zariski open subset $U$
of $X$, $Z\cap (U\times Y)$ is the graph of a holomorphic
function from $U$ to $Y$.

Fujiki, in his study \cite{Fujiki} of automorphism groups of
compact K\"{a}hler manifolds, introduces the notion of a
``meromorphic group". As we will be proposing a less
restrictive meaning for ``meromorphic group" we will refer
to Fujiki's notion as ``Fujiki-meromorphic". A
Fujiki-meromorphic group is a complex Lie
group $G$ which is a Zariski open subset of some compact
complex space $G^{*}$ such that the group operation of
$G$ extends to a meromorphic mapping from $G^{*}\times
G^{*}$ to $G^{*}$ which is holomorphic on $(G\times
G^{*})\cup(G^{*}\times G)$. Let ${\mathcal C}$ be the full
subcategory of
${\mathcal A}$ consisting of compact complex spaces which are
holomorphic images of compact  K\"{a}hler
manifolds. Fujiki proves that if $G$ is a Fujiki-meromorphic
group in $\mathcal C$ (namely $G^{*}\in {\mathcal C}$), then $G$ is
``meromorphically" isomorphic to an extension of a complex
torus by a linear algebraic group, generalizing Chevalley's
well-known theorem for algebraic groups. He raises the
issue whether this remains true in the more general
category ${\mathcal A}$ and proves it for $G$ commutative.

Our proposed definition of a meromorphic group is
as follows: $G$ is a connected complex Lie group with a
finite covering by Zariski open subsets $U_{i}$ of
irreducible compact complex spaces
$X_{i}$  ($i=1,..,n$) such that both the transition maps and
the group operation on $G$ extend to meromorphic maps
between the various $X_{i}$ and their products. Note that
if the $X_{i}$ happen to be algebraic varieties then this
agrees with the definition of an abstract algebraic group.
Complex
algebraic groups, complex tori, and Fujiki-meromorphic
groups are all meromorphic groups. In fact our
results imply that meromorphic groups coincide with
Fujiki-meromorphic groups, and moreover have K\"{a}hler
compactifications.

This paper is informed by model-theoretic concerns, and
indeed some model-theoretic results play a role in the
proofs. The main point is that ${\mathcal A}$ considered as a
many-sorted structure, whose sorts are the compact complex
spaces and whose basic relations are the analytic subsets of
various Cartesian products of the sorts, is a structure with
quantifier-elimination and finite Morley rank (sort by sort).
This was proved by Zilber \cite{Zilber} (although
quantifier-elimination was also noted earlier in
\cite{Loj}). Quantifier-elimination says that the definable
sets are precisely the finite unions of locally Zariski
closed subsets of various compact complex spaces. It follows
that definable functions are precisely ``piecewise
meromorphic" functions. Moreover definable groups are (up to
definable isomorphism) precisely the meromorphic groups,
giving another justification for introducing the notion of
meromorphic group. A strongly minimal set in ${\mathcal A}$
  is a definable set without infinite, co-infinite definable
subsets. A strongly minimal group in ${\mathcal A}$ is precisely a
meromorphic group without proper infinite Zariski-closed
subsets. In  \cite{Hrushovski-Zilber1} it was noted that the
deep results of \cite{Hrushovski-Zilber2} apply to strongly
minimal sets definable in ${\mathcal A}$, implying that any
strongly minimal definable group $G$ is either an
($1$-dimensional) algebraic group, or is ``modular": every
definable subset of $G \times \cdots \times G$ is essentially a
translate of a subgroup. Simple complex tori of dimension $>
1$ are examples of strongly minimal modular groups (see
\cite{Pillay}). For the converse, Hrushovski
\cite{Hrushovski} asked whether strongly minimal
modular groups are (necessarily simple) complex tori. In fact
Hrushovski outlined to the first author some ideas for
proving this, depending however on finding a good
compactification of the group. In any case, the question was
answered by the second author in \cite{Scanlon} for the
special case when $G$ is itself interpretable in a strongly
minimal compact complex manifold. In
\cite{Pillay} additional observations about ${\mathcal A}$ and
its model theory were made, including ``elimination of
imaginaries". Also, it was asked whether the Chevalley
theorem holds for groups definable in ${\mathcal A}$. We found
subsequently that the same question was asked in
\cite{Fujiki} for Fujiki-meromorphic groups.

We will prove the following results:
\begin{Theorem} Suppose $G$ is a strongly minimal
meromorphic group. Then $G$ is
meromorphically isomorphic to either a $1$-dimensional
algebraic group or a simple nonalgebraic torus.
\end{Theorem}

\begin{Theorem} Suppose $G$ is a connected meromorphic group.
Then $G$ has a normal connected meromorphic subgroup $L$
such that $L$ is (meromorphically isomorphic to) a
linear algebraic group, and $G/L$ is a complex torus.
Moreover $L$ is unique.
\end{Theorem}

Theorem 1.1 is a special case of Theorem 1.2. Theorem 1.1
will be proved by finding a good compactification of
$G$ (i.e. showing that
$G$ is Fujiki-meromorphic) and then (as $G$ is commutative)
referring to \cite{Fujiki}. By again finding a suitable
compactification we will first prove Theorem 1.2 for the
special case when $G$ is an extension of a $1$-dimensional
linear algebraic group by a simple complex torus. The general
case will follow by an induction on dimension, making use of
some additional ingredients such as the structure of
compact complex spaces with algebraic co-dimension $1$, and
some model theory of groups of finite Morley rank.

In the next section we give some definitions and recall both
complex analytic and model-theoretic notions. In section 3 we
carry out compactifications, proving Theorem 1.1 among
other things. In section 4 we prove Theorem 1.2. Some
additional remarks are made in section 5.

\section{Preliminaries}
For basic results, notions and notation concerning complex
spaces and meromorphic maps, we refer the reader to
\cite{Fischer}, \cite{Grauert-Peternell-Remmert} and
\cite{Ueno}. However we will repeat a few crucial
definitions and results which we will be relying on. For us
${\mathcal A}$ denotes the class of reduced irreducible compact
complex spaces. We take as given the notion of a holomorphic
map $f$ from $X$ to $Y$ where
$X,Y \in {\mathcal A}$. $\dim(X)$ denotes the complex dimension of
$X$. A {\em modification} of
$X\in {\mathcal A}$ is some $Y\in {\mathcal A}$ and a surjective
holomorphic
$f:Y\rightarrow X$ such that for some proper closed analytic
subsets $A$ of $Y$ and $B$ of $X$,
$f|(Y\setminus A):{Y\setminus A}\rightarrow X\setminus B$
is biholomorphic. Resolution of singularities says that any
$X$ has a modification $(Y,f)$ such that $Y$ is nonsingular
(so a connected compact complex manifold). The notion of a
meromorphic mapping $f$ from $X$ to $Y$ ($X,Y\in {\mathcal A}$)
is crucial. Such an object can be defined in various
equivalent ways. For $X$ irreducible we define $f$ to be a
function from $X$ to the set of subsets of $Y$ such that (i)
the ``graph" of $f$, $\{(x,y)\in X\times Y: y\in f(x)\}$ is
an irreducible analytic subset of $X\times Y$, and for all
$x$ in some (dense) Zariski open subset $U$ of $X$, $f(x)$ is
a singleton. For a general $X$ we say that $f$ is meromorphic
if each of its restrictions to the irreducible components of
$X$ is meromorphic. We say that $f$ is holomorphic, or
defined, at the points in $U$. Let $Z$ be the graph of $f$
as defined above, and $\pi$ the projection from $Z$ onto
$X$. Then $(Z,\pi)$ turns out to be a modification of $X$.
The projection of $Z$ on the second coordinate is then a
holomorphic map from $Z$ to $Y$ which is said to be
"resolution of indeterminacies" of the meromorphic map $f$.
 From the definition of a meromorphic map one easily
derives the following fact.

\begin{Fact} Let $X,Y\in {\mathcal A}$. Let $f,g$ be meromorphic
mappings from $X$ to $Y$. Suppose that for some dense
Zariski open subset $U$ of $X$, $f$ and $g$ agree on $U$.
Then $f=g$.
\end{Fact}

Suppose that $U$ is a dense Zariski open subset of $X\in
{\mathcal A}$, and $f$ a holomorphic map from $U$ into
$Y\in{\mathcal A}$. By abuse of language we may sometimes say
that $f$ is meromorphic if there is a meromorphic mapping
$g$ from $X$ to $Y$ which agrees with
$f$ on $U$. A natural category which can be associated to
${\mathcal A}$ is the category whose objects are those complex
spaces which are Zariski open subsets of spaces in ${\mathcal
A}$ and whose morphisms are the holomorphic maps which are
meromorphic in the sense of the previous sentence. If we
restrict our attention to those $X$, $Y$ which are
projective algebraic varieties, this category is exactly
that of quasiprojective varieties and morphisms.

It is also natural to consider complex spaces which have a
finite covering by Zariski open subsets $U_{i}$ of spaces
$X_{i}$ in
${\mathcal A}$ where the transition maps are meromorphic in the
above sense. Morphisms in this category would be holomorphic
maps which are meromorphic (in the above sense) when read in
each $U_{i}$. What we will call a {\em meromorphic} group
is exactly a group object in this latter category. Here is
the precise definition.

\begin{Definition} A meromorphic group is a
connected complex Lie group $G$, with a finite covering by
open subsets $W_{i}$, for $i=1,..,n$, and for each $i$ a
(biholomorphic) isomorphism $\phi_{i}$ of $W_{i}$ with a
Zariski open subset $U_{i}$ of some $X_{i}\in {\mathcal A}$ such
that
\begin{itemize}
\item[(i)] For each $i\neq j$, $\phi_{i}(W_{i}\cap W_{j})$ is a
Zariski open subset of $X_{i}$, and the induced biholomorphic
map between $\phi_{i}(W_{i}\cap W_{j})$ and
$\phi_{j}(W_{i}\cap W_{j})$ is meromorphic, namely is the
restriction of a meromorphic mapping between $X_{i}$ and
$X_{j}$.
\item[(ii)] For each $i,j,k$, $\{(x,y)\in U_{i}\times
U_{j}:\phi_{i}^{-1}(x)\cdot\phi_{j}^{-1}(y) \in W_{k}\}$ is
Zariski open in $X_{i}\times X_{j}$ and the induced
holomorphic map ((x,y) goes to $\phi_{k}(\phi_{i}^{1}(x)\cdot
\phi_{j}^{-1}(y))$) from
$U_{i}\times U_{j}$ to $U_{k}$ is meromorphic, namely is the
restriction of a meromorphic mapping between $X_{i}\times
X_{j}$ and $X_{k}$.
\end{itemize}
\end{Definition}

Conditions (i) and (ii) can be expressed briefly by saying
that the transition maps as well as the group operation are
meromorphic when read in the various $U_{i}$ and their
Cartesian products.

We  say that the covering by the
$W_{i}$'s and the isomorphisms with the $U_{i}$'s satisfying
(i) and (ii) above give the complex Lie group $G$
a meromorphic structure.

If $G$ is a meromorphic group as in Definition 2.2,
by a meromorphic subgroup $H$ of $G$ we mean a closed
subgroup such that for each $i$, $\phi_{i}(H\cap W_{i})$ is
the intersection of an analytic subset of $X_{i}$ with
$U_{i}$. Clearly $H$ has the structure of a meromorphic
group.

A holomorphic homomorphism (complex Lie homomorphism) $f$
between meromorphic groups
$G_{1}$ and $G_{2}$ is meromorphic if when restricted to
the charts the map is meromorphic, that is, extends to
meromorphic mappings between the relevant compact complex
spaces.

So now we have the category of meromorphic groups and
meromorphic homomorphisms. The following says that quotient
objects exist. It follows from looking at the equivalent
category of definable groups and using the elimination of
imaginaries result from \cite{Pillay}. This will be explained
below.

\begin{Fact} Let $G$ be a meromorphic group and $N$ a normal
meromorphic subgroup. Then there is a meromorphic group $H$
and a surjective meromorphic homomorphism from $G$ to $H$
whose kernel is $N$.
\end{Fact}

We now repeat the definition of a Fujiki-meromorphic group
and recall what Fujiki proved.
\begin{Definition} Let $G$ be a complex Lie group.
\begin{itemize}
\item[(i)] A meromorphic compactification of $G$ is a compact
complex space
$G^{*} \in {\mathcal A}$ which contains $G$ as a dense Zariski
open subset, such that the group operation $\mu:G\times G
\rightarrow G$ is meromorphic, i.e. the restriction of a
meromorphic mapping $\mu^{*}$ say, from $G^{*}\times G^{*}$
to $G^{*}$

\item[(ii)] A Fujiki-compactification of $G$ is a
meromorphic compactification
$(G^{*},\mu^{*})$ of $G$ such that $\mu^{*}$ is holomorphic
on $(G\times G^{*})\cup (G^{*}\cup G)$.

\item[(iii)] $G$ is Fujiki-meromorphic if $G$ has a
Fujiki-compactification.
\end{itemize}
\end{Definition}

\begin{Remark} (i) A Fujiki-meromorphic group is a
meromorphic group.
\newline
(ii) A connected compact complex Lie group (i.e. a complex
torus) is Fujiki-meromorphic.
\newline
(iii)(Remark 2.3 of \cite {Fujiki}.) A complex
algebraic group is Fujiki-meromorphic.
\end{Remark}

Following notation of Fujiki \cite{Fujiki}:
\begin{Definition} We will call the meromorphic
group
$G$ regular if there is a meromorphic homomorphism $f$ from
$G^0$, the connected component of the identity in $G$,
onto a complex torus $T$ such that the kernel $L$ of $f$
is meromorphically isomorphic to a connected linear algebraic
group. (Briefly said: $G^0$ is meromorphically an extension
of a complex torus $T$ by a linear algebraic group $L$.)
\end{Definition}

\begin{Remark}
\begin{itemize}
\item[(i)] Let $G$ be regular and let $T, L$ be as
above. Then $L$ and $T$ are unique. In particular $L$ is the
unique maximal normal connected meromorphic subgroup of $G^0$
which is meromorphically isomorphic to a linear algebraic
group.

\item[(ii)] A regular meromorphic group is
Fujiki-meromorphic.
\end{itemize}
\end{Remark}
\begin{proof}
  (i) Suppose $L_{1}$ is a normal connected
meromorphic subgroup of $G$ which is meromorphically
isomorphic to a linear algebraic group. Then $L_{1}/L$
meromorphically embeds in $T$. So $L_{1}/L$ is both a
complex torus and a linear algebraic group forcing it to be
trivial. That is $L_{1}$ is contained in $L$.

(ii) As in Remark 2.3 of \cite{Fujiki}.
\end{proof}

Recall that ${\mathcal C}$ is the subclass (in fact full
subcategory) of ${\mathcal A}$ consisting of those $X$ which are
holomorphic images of compact connected K\"{a}hler manifolds.
We will say that the connected
meromorphic group $G$ is of type ${\mathcal C}$ if there is a
Fujiki-compactification $G^{*}$ of $G$ which is in ${\mathcal
C}$. Fujiki proves:
\begin{Fact} (i) Suppose $G$ is
Fujiki-meromorphic. Then $G$ is regular iff $G$ is of type
${\mathcal C}$.
\newline
(ii) Suppose that $G$ is commutative and
Fujiki-meromorphic. Then $G$ is regular.
\end{Fact}

In the final part of this section we discuss the model
theory of compact complex manifolds. We will have to assume
the basics of model theory, and a bit more. \cite{Hodges} is
a good reference for basic model theory. The first four
chapters of \cite{Bouscaren} (by Bouscaren, Ziegler,
Lascar, Pillay) are a useful reference for various aspects of
applied and geometric stability theory.
\cite{Pillay-book} is an advanced text on geometric
stability. \cite{Borovik-Nesin} deals with the theory of
groups of finite Morley rank. Another good reference for
stable groups is \cite{Poizat}.

We consider
${\mathcal A}$ as a many-sorted first order structure whose
sorts are the (reduced, irreducible) compact complex spaces
and basic relations the analytic subsets of finite Cartesian
products of such things.
\begin{Fact} $Th({\mathcal A})$ has quantifier-elimination,
elimination of imaginaries and each sort has finite Morley
rank. Moreover ${\mathcal A}$ is $\aleph_1$-compact.
\end{Fact}
Quantifier-elimination was proved in \cite{Loj}, and
independently in \cite{Zilber}. It says that any
definable subset of a sort $X$ is ``analytically
constructible", that is a finite union of intersections of
analytic (Zariski closed) subsets and complements of
analytic (Zariski open) sets. A characterization of
definable functions follows from this: Suppose $U$ is a
definable set, and $f$ a definable function from $X$ into
some sort $Y$. Then we can write $X$ as a disjoint union of
definable sets $U_{i}$, where each $U_{i}$ is a Zariski open
subset of some sort (complex space) $X_{i}$ such that for
each
$i$, the restriction of
$f$ to $U_{i}$ is holomorphic and is the restriction to
$U_{i}$ of a meromorphic mapping from $X_{i}$ to $Y$. We say,
with possibly some abuse of language, that definable
functions are piecewise meromorphic.

Zilber \cite{Zilber} proved finiteness of Morley rank.
Elimination of imaginaries was observed in
\cite{Pillay}.

$\aleph_{1}$-compactness means that any countable family of
definable subsets of some sort $X$ has nonempty intersection
as long a every finite subfamily does.

With Fact 2.9 there is a remarkable parallel
between complex-analytic and model-theoretic structural
and classification results. We refer the reader to
\cite{Zilber}, \cite{Hrushovski}, \cite{Pillay} for
more discussion. Note that on the face of it ${\mathcal A}$ is
not $\aleph_{1}$-saturated, as each element of each sort is
essentially named by a constant. One can ask whether there
is some sublanguage $L_{0}$ of the full language $L$
described above such that every relation in $L$ is definable
possibly with parameters in the language $L_{0}$ and such
that the reduct ${\mathcal A}|L_{0}$ {\em is}
$\aleph_{1}$-saturated. This is not true, as, for example,
a general generalized Hopf surface has continuum many
holomorphic automorphisms but our Proposition 5.2 shows that
it has trivial generic type and hence
cannot have a non-trivial definable family of
automorphisms. On the other hand,
${\mathcal C}$ can be considered as a reduct of ${\mathcal A}$
(fewer sorts but the full structure on each sort), which
has quantifier-elimination and elimination of imaginaries
in its own right. Fujiki's results \cite{Fujiki2} on the
Douady spaces of manifolds in ${\mathcal C}$ imply that the
structure ${\mathcal C}$ {\em is} $\aleph_{1}$-saturated in a
suitable sublanguage. (See \cite{Moosa}.)

${\mathcal A}'$ will denote a very saturated elementary
extension of ${\mathcal A}$. For any $X\in {\mathcal A}$, $X'$
denotes its canonical extension in ${\mathcal A}'$. We will
often work model-theoretically in ${\mathcal A}'$. For example,
a definable property holds generically on $X$ iff it holds
of a generic point of $X'$.

A definable group
$G$  in
${\mathcal A}$ will be called connected if
$G$ has no definable subgroups of finite index. Any
meromorphic group is clearly a definable group (using
elimination of imaginaries). Methods from the algebraic case
due to Hrushovski and van den Dries (see
\cite{Poizat} as well as Pillay's article in
\cite{Bouscaren})
   adapt to yield the important:
\begin{Fact} Any  definable connected group $G$ in ${\mathcal
A}$ is definably isomorphic to a connected meromorphic
group $H$ (unique up to meromorphic isomorphism).
\end{Fact}
This fact gives a natural equivalence between the category
of definable groups and meromorphic groups. In particular
any definable homomorphism between meromorphic groups will
be meromorphic (in particular holomorphic). From here on we
will use ``definable" interchangeably with ``meromorphic"
when talking about groups and homomorphisms.

A definable set $X$ (in ${\mathcal A}$) is said to be strongly
minimal if $X$ is
infinite and has no infinite co\"{\i}nfinite definable subsets. A
definable connected group
$A$ is said to be {\em modular} if every definable subset of
$A^{n}$ is a Boolean combination of translates of definable
subgroups. In
\cite{Hrushovski-Zilber1} it was proved that the results of
\cite{Hrushovski-Zilber2} apply to the category
${\mathcal A}$. This yields.
\begin{Fact} Suppose $G$ is a definable connected group in
${\mathcal A}$ which has no infinite normal definable subgroups.
Then either $G$ is strongly minimal and modular or $G$ is
definably isomorphic to a (complex) algebraic group.
\end{Fact}

It follows that if $T$ is a nonalgebraic simple complex
torus then $T$ is modular. (A direct proof, avoiding
\cite{Hrushovski-Zilber2} was given in \cite{Pillay}.)

\section{Compactifications}
We will prove:
\begin{Theorem} Let $G$ be a connected commutative
meromorphic group, which is either strongly minimal, or an
extension of a connected $1$-dimensional linear algebraic
group by a simple complex torus. Then $G$ is
Fujiki-meromorphic.
\end{Theorem}

A consequence is:
\begin{Corollary}
\begin{itemize}
\item[(i)] Let $G$ be a strongly minimal
meromorphic group. Then $G$ is meromorphically isomorphic
to either a $1$-dimensional algebraic group or a
simple modular complex torus.
\item[(ii)] Let $G$ be a commutative meromorphic group which is
an extension of a $1$-dimensional linear algebraic group by
a simple complex torus. Then $G$ meromorphically splits.
\end{itemize}
\end{Corollary}
\begin{proof}
(i) $G$ is commutative, so by Theorem 3.1, Fujiki
meromorphic, thus by Fact 2.8 (ii), meromorphically an
extension of a complex torus $T$ by a linear algebraic
group $L$. As $G$ is strongly minimal, $G$ is either $T$ or
$L$. If $G = L$, then $\dim(L) = 1$. If $G = T$, then $T$ has
no proper infinite analytic subsets so is either an elliptic
curve or  simple and modular (by 2.11).

(ii) Immediate, by Fact 2.8 (ii).
\end{proof}

To prove Theorem 3.1 we will find a
meromorphic compactification
$G^{*}$ of $G$ and then show it to be a
Fujiki-compactification. The following general result
concerning compactifications of commutative meromorphic
groups will be useful.

\begin{Lemma} Suppose that the connected commutative
meromorphic group
$(G,\mu)$ has meromorphic compactification $(G^{*},\mu^{*})$.
   Suppose $S = G^{*}\setminus G$ is nonempty. Then
\begin{itemize}
\item[(i)] Every component of $S$
has co-dimension $1$ in $G^{*}$.
\item[(ii)] $\mu^{*}|(G^{*}\times S)$ is a meromorphic mapping from
$G^{*}\times S$ to $S$.
\item[(iii)] For each $g\in G$, and component $C$ of $S$,
$\mu^{*}_{g} = \mu^{*}(g,-):C\rightarrow C$ is biholomorphic
on a dense Zariski open subset of $C$, and for $g,h\in G$,
$\mu^{*}_{g}.\mu^{*}_{h} = \mu^{*}_{g.h}$ on a dense
Zariski open subset of $C$.
\end{itemize}
\end{Lemma}
\begin{proof}
(i) Let $n =
\dim(G^{*})$ (=$\dim(G)$). Suppose for the sake of
contradiction that there is $x\in S$ such that
$\dim_{x}(S) < n-1$. Let $\Delta _{n}$ be the
open unit disc in ${\mathbb C}^{n}$ and $f: \Delta _{n}
\rightarrow U$ be a coordinate function for any open
neighborhood $U$ of $x$ in $G^{*}$ where $U$ is chosen such
that $U\cap S$ has dimension $< n-1$. So if $A =
f^{-1}(U\cap S)$, then $A$ is an analytic subset of $\Delta
_{n}$ of codimension at least $2$, and $f_{1}$ = $f|(\Delta
_{n}\setminus A)$ is a holomorphic embedding into $G$. As $G$
is a connected commutative Lie group its universal cover is
$\pi:{\mathbb C}^{n}\rightarrow G$. As $A$ has co-dimension at
least two in $\Delta _{n}$, $\Delta _{n}\setminus A$ is
simply connected, and so $f_{1}$ lifts to a holomorphic map
$f_{2}: \Delta_{n} \rightarrow {\mathbb C}^{n}$ (see
\cite{Fischer}). Let
$g = \pi \circ f_{2}$. Then $g$ is a holomorphic map from $\Delta$
into $G^{*}$ which agrees with $f$ off the thin analytic
subset $A$. But then by Riemann's removable singularity
theorem (\cite{Fischer})
$f = g$, contradicting the fact that
$x\notin G$. (i) is proved.

(ii) Let $\Gamma$ be the graph of $\mu^{*}$. We will first
show that for  all
$(g,x)$ in some dense Zariski open subset
$V$ of $G^{*}\times S$, $\{y: (g,x,y)\in \Gamma\}$ is finite.
If not, then for a Zariski open subset $V$ of $G^{*}\times
S$, the above set of $y$'s has positive dimension. It
follows from (i) that $\dim(\Gamma \cap (G^{*}\times S\times
G^{*}) \geq 2n$, contradicting irreducibility of $\Gamma$.
It now follows by the implicit function theorem that

(a) $\mu^{*}$ is holomorphic on $V$.

For $g\in G$ let $\mu^{*}_{g}$ be $\mu^{*}(g,-)$, a
meromorphic mapping from $G^{*}$ to $G^{*}$. Note that if
$\mu^{*}_{g}$ is defined (single valued) at $x$ and
$\mu^{*}_{h}$ is defined at $\mu^{*}_{g}(x)$ then
$\mu^{*}_{hg}$ is defined at $x$ and equals
$\mu^{*}_{h}(\mu^{*}_{g}(x))$. It follows that if $(g,x)\in
V$ and $g\in G$, then $\mu(g,x)\in S$. As $G$ is
Zariski-dense in $G^{*}$ it follows that

(b) $\mu^{*}|(G^{*}\times S)$ is a meromorphic mapping into
$S$, yielding (ii).

(iii) The same argument as above shows that for any $g\in G$,
$\mu^{*}_{g}|S$ is a meromorphic mapping from $S$ to $S$.
Let $C_{1},..,C_{s}$ be the irreducible components of $S$.
Note that the image of the meromorphic mapping $\mu^{*}_{g}$
from $G^{*}$ to $G^{*}$ (i.e. projection of its graph on
second component) is all of $G^{*}$. But for $x\in G$,
$\mu^{*}_g(x)\in G$. Thus the image of the meromorphic mapping
$\mu^{*}_{g}|S$ is all of $S$. We work model-theoretically.
Fix $C_{i}$. Let $x$ be a generic point of $C_{i}'$ over
${\mathcal A}$. So $y= \mu^{*}_{g}(x)\in C_{j}'$ for some $j=
f(i)$. It follows that $\mu^{*}_g|C_{i}$ is a meromorphic
mapping from $C_{i}$ into $C_{f(i)}$.

(c) Thus $f = f_{g}$ must be a permutation of $\{1, \ldots, s \}$.

If for some $i$, and $x$ as above, $\mu^{*}_{g}(x)$ is not a
generic point of $C_{f(i)}'$ over ${\mathcal A}$, then there is a
proper analytic subset $D_{f(i)}$ of $C_{f(i)}$ such that
$\mu^{*}_{g}|C_{i}$ has image contained in $D_{f(i)}$. By
(c), we contradict the fact that $\mu^{*}_{g}|S$ has image
all of $S$.

Thus for $x\in C_{i}'$ generic, $\mu^{*}_{g}(x)$ is generic
in $C_{f(i)}'$ over ${\mathcal A}$. It follows that $g\rightarrow
f_{g}$ gives a definable action of $G$ on $\{1, \ldots, s\}$. As
$G$ is connected this has to be trivial. This gives (iii).
\end{proof}

\begin{Remark}
(iii) above is interpreted
model-theoretically by saying that $G$ acts generically on
$C$: let $p = p_{C}$ be the generic type of the component
$C$ of $S$. Then for $g\in G'$ and $x$ realizing $p$
independent of $g$ (over ${\mathcal A}$), $\mu^{*}(g,x)$ is
defined, realizes $p$ and is independent from $g$. Moreover,
if $g,h\in G'$ and $x$ realizes $p$ independent of $g,h$
then $\mu^{*}(hg,x) = \mu^{*}(h,\mu^{*}(g,x))$.
\end{Remark}

We can now obtain the strongly minimal case of Theorem 3.1.
\begin{Lemma} Let $G$ be a strongly minimal meromorphic
group. Then $G$ is Fujiki-meromorphic.
\end{Lemma}
\begin{proof}
\newline
\emph{Step 1.} Finding a meromorphic compactification.

By assumption on
$G$ some open nonempty definable subset $U$ of $G$ is
already a Zariski open subset of a compact complex space
$X$, which we may assume by resolution of singularities to be
a manifold. Moreover by strong minimality of
$G$, $G\setminus U$ is finite, say $\{g_{1}, \ldots, g_{n}\}$. For
$i=1, \ldots, n$ let $V_{i}$ be a coordinate neighborhood of
$g_{i}$ in $G$ such that the closures ${\bar V_{i}}$ of the
$V_{i}$ in $G$ are disjoint. Note that $ K_{i} = {\bar
V_{i}}\setminus \{g_{i}\}$ is contained in $U$, so in $X$,
but is not compact, so not closed in $X$. Let $D_{i}$ be
the boundary of $K_{i}$ in $X$, namely ${\bar
K_{i}}\setminus K_{i}$. Let $\pi:X\rightarrow X'$ be the
quotient map which collapses each $D_{i}$ to a point $c_{i}$.
Then $X'$ is compact, $\pi$ is holomorphic (in fact is a
modification), and is biholomorphic outside the union of the
$D_{i}$'s. Let $f:G\rightarrow X'$ be defined by $f(x) =
\pi(x)$ for $x\in U$ and $f(g_{i}) = c_{i}$. Then $f$ is a
definable, holomorphic embedding.

\noindent
\emph{Step 2.} $G^{*}$ is a Fujiki-compactification of $G$.

If $G = G^{*}$
there is nothing to do. Otherwise, (as $G$ is commutative)
Lemma 3.3 applies. Let $S$ be as there. We will show that
$G$ is holomorphic on
$S$, and in fact acts as the identity. Note that the generic
type of $G$ is orthogonal to any set of dimension less than
that of $G$ ($G$ being strongly minimal). In particular $G$
is orthogonal to $S$. Fix a component $C$ of $S$. Lemma 3.3
(iii) gives us a generic action of $G$ on $C$. Let
$g,h\in G'$ be generic independent elements of $G$ and let
$x$ be generic in $C'$ over $\{g,h\}$. Then by the
orthogonality mentioned above, each of $g$ and $h$ is
independent from
$\{x,\mu^{*}(g,x)\}$. It follows that $\mu^{*}(g,x) =
\mu^{*}(h,x)$, and thus $\mu^{*}(h^{-1}.g,x) = x$. But
$h^{-1}.g$ is generic in $G'$ and independent from $x$. It
follows that $G$ acts generically trivially on $C$. So the
holomorphic map from $G^{*}\times C$ to $C$ taking $(g,x)$
to $x$ agrees generically with the meromorphic mapping
$\mu^{*}|(G^{*}\times C): G^{*}\times C\rightarrow C$. By
Fact 2.1, these mappings agree. This shows that
$(G^{*},\mu^{*})$ is a Fujiki-compactification of $G$.
\end{proof}

We now deal with the case when $G$ is a commutative
extension of the additive group, ${\mathbb G}_{a}$,
or  the multiplicative group, ${\mathbb G}_{m}$, by a simple complex torus
$T$.  We let
$H$ denote $G/T$ (so $H$ is ${\mathbb G}_a$ or ${\mathbb G}_m$). If
$G$ is meromorphically isomorphic to an algebraic group,
then $G$ is clearly Fujiki-meromorphic (in fact the
Chevalley theorem applies immediately, yielding Theorem
1.2). If $T$ has a definable complement in $G$ (up to
finite), then again we get the required conclusion. So  for
the rest of this section we make:

\noindent
\emph{Assumption.}
\begin{itemize}
\item[(a)] $G$ is a commutative meromorphic
group, which is meromorphically an extension of $H$ by a
simple complex torus $T$, where $H$ is ${\mathbb G}_a$ or ${\mathbb G}_m$.

\item[(b)] $G$ is not meromorphically isomorphic to an algebraic
group.

\item[(c)] There is no definable connected infinite subgroup $L$
of $G$ with $L\cap T$ finite.
\end{itemize}
\bigbreak

We will show that $G$ is Fujiki-meromorphic (which actually
leads to a contradiction).

We will make use of the ``socle theory" from
\cite{Hrushovski-ML}.

\begin{Lemma} $T$ is the maximal almost pluriminimal
definable subgroup of $G$.
\end{Lemma}
\begin{proof} Note that $T$ being simple is almost strongly
minimal. So if the lemma, as $G/T$ has
dimension $1$, $G$ is semipluriminimal. By
\cite{Hrushovski-ML}, $G$ is an almost direct product of
pairwise orthogonal semiminimal groups. If $G$ is already
semiminimal, then as $G$ is nonorthogonal to ${\mathbb P}^{1}$ via
$G\rightarrow H$, $G$ must be algebraic, contradicting
Assumption (b) Thus $G$ is the semidirect product of $T$ and
some $L$, contradicting Assumption (c).
\end{proof}

\begin{Lemma} Let $X$ be a definable subset of $G$. Assume
that the Morley rank of $X$ is strictly less than the Morley
rank of $G$ (equivalently $X$ is not Zariski-dense in $G$).
Then $X$ is contained in finitely many translates of $T$.
\end{Lemma}
\begin{proof} We prove the lemma by induction on $RM(X) = m$.
It is clearly true for $m=0$. We may assume that the Morley
degree of $X$ is $1$. Let $S$ be the (model-theoretic)
stabilizer of $X$. If $S$ is finite then by Lemma 3.6, as
well as Proposition 4.3 of \cite{Hrushovski-ML}, $X$ is, up
to a set of Morley rank $< m$, contained in a single
translate of $T$. By induction hypothesis, $X$ is contained
in finitely many translates of $T$ as desired. So we may
assume that $S$ is infinite. By Assumption (c), and the fact
that $T$ is simple, $S$ must contain $T$. Note that $RM(T) =
RM(G) - 1 \geq RM(X)$. But it is well-known that the Morley
rank of the stabilizer of a Morley degree $1$ set $X$ is at
most the Morley rank of $X$, with equality if and only if the
stabilizer is connected and $X$ is, up to a set of smaller
Morley rank, a translate of this stabilizer. Thus
$RM(S) = RM(X)$, $S = T$ and up to a set of
smaller Morley rank, $X$ is a translate of $T$, so we finish
again by induction.
\end{proof}

\begin{Lemma} $G$ is Fujiki-meromorphic.
\end{Lemma}
\begin{proof} As in the strongly minimal case we first find a
compact complex manifold $G^{*}$ containing $G$ as a Zariski
open set, and then show that this gives $G$ a
Fujiki-meromorphic structure.

\noindent
\emph{ Step I.} Finding the compactification.

Let $RM(G) = n$.  By definition of $G$ being
a meromorphic group, let $U$ be a definable subset (with
Morley rank $n$) of $G$ which is a dense
Zariski-open subset of a compact complex manifold ${\bar U}$.
Let $\pi: G \rightarrow H$ be the canonical surjective
homomorphism. Then $\pi$ takes $U$ onto a cofinite subset
$\pi(U)$ of $H$.

\noindent
\emph{Claim 1.} We may assume that for any $x\in \pi(U)$,
$\pi^{-1}(x)\cap U$ = $\pi^{-1}(x)$ (a translate of $T$).

\begin{proof}
$Y$ = $\pi^{-1}(\pi(U))\setminus U$ is a
definable subset of $G$ of Morley rank $< n = RM(G)$. By
Lemma 3.7, $Y$ is contained in finitely many translates of
$T$, namely finitely many fibers of $\pi$. Remove these
fibers from $U$.
\end{proof}

Let $\pi'$ denote $\pi|U$. $\pi'$ extends to a meromorphic
function ${\bar \pi}$ from ${\bar U}$ to ${\mathbb P}^{1}$.
Further restricting $U$ we may assume:

\noindent
\emph{Claim 2.} For all $x\in \pi'(U)$, $(\pi')^{-1}(x) =
{\bar
\pi}^{-1}(x)$.
\bigbreak

Let $C$ be the finite set $H\setminus \pi(U)$. Then we can
find $h\in\pi(U)$ such that $h.C\subset \pi(U)$. Let $g\in
U$ be a preimage of $h$. Let $\tau_{g}:G\rightarrow G$ be
multiplication by $g$. $\tau_{g}|U$ is not defined
everywhere but is holomorphic on the open set where it is
defined and so extends to a meromorphic map ${\bar \tau_{g}}:
{\bar U}\rightarrow {\bar U}$. By a theorem of Remmert (see
Theorem 1.9 in Chapter VII of
\cite{Grauert-Peternell-Remmert}), there is a modification
$\nu: {\tilde U}\rightarrow {\bar U}$ and a holomorphic map
${\tilde \tau}: {\tilde U}\rightarrow {\bar U}$ such that
${\tilde \tau} = {\bar \tau_{g}}\cdot \nu$. In particular,
for $x\in {\tilde U}$ such that $\tau_{g}|U$ is defined at
$\nu (x)$,
${\tilde \tau} (x) = \tau_{g}(\nu (x))$.

\bigbreak
\noindent
\emph{Claim 3.} ${\bar \pi}\cdot {\tilde \tau} =
\tau_{h}\cdot {\bar \pi}\cdot \nu$.
\begin{proof}
  This holds generically, so holds everywhere.
\end{proof}

We will now construct the required compactification $G^{*}$
of $G$ as a holomorphic image of ${\tilde U}$. Let $S =
{\mathbb P}^{1}\setminus H$. So $S = \{\infty\}$ or
$\{\infty, 0\}$. As a set $G^{*}$ will be the disjoint union
of $G$ with
${\bar \pi}^{-1}(S)$. The manifold structure of $G^{*}$ is
as follows: $G$ is given its canonical manifold structure.
Now let $x\in {\bar \pi}^{1}(S)$. Let $y = {\bar
\pi}(x) \in {\mathbb P}^{1}$. Choose an open neighborhood $V$ of
$y$ in ${\mathbb P}^{1}$ such that $V\setminus\{y\} \subset U$.
Then ${\bar \pi}^{-1}(V) \subset G^{*}$ is an open
neighborhood of $x$. The transition maps are clearly
holomorphic, yielding a structure of a complex compact
manifold on $G^{*}$ containing $G$ as an open (dense)
subset.

Now we define a holomorphic surjective map $f$ from ${\tilde
U}$ to $G^{*}$. Let $x\in {\tilde U}$. If
${\bar \pi}({\tilde \tau}(x)) \notin C$ define $f(x) =
{\tilde \tau}(x)$ (so $f(x) \in {\bar \pi}^{-1}(\pi(U)\cup
S)\subset G^{*}$). On the other hand, if
${\bar \pi}({\tilde
\tau}(x)) \in C$, define $f(x) = g.\nu(x)$. Note that in
this latter case $\nu(x) \in U\subset G$ and so $g.\nu(x)
\in G$ and $\pi(g.\nu(x)) = {\bar \pi}{\tilde \tau}(x)$.

It is easy to check, given our assumptions, that $f$ is
holomorphic and surjective. So $G^{*}$ is a compact complex
manifold, containing $G$ as a dense Zariski open subset (the
embedding of $G$ in $G^{*}$ is definable and holomorphic).

This completes Step I.

\bigbreak
\noindent
\emph{Step II:} The action of $G$ on itself
extends to a trivial action on the boundary
$G^{*}\setminus G$.

Let $C_{1}, \ldots, C_{k}$ be the irreducible components of
$G^{*}\setminus G$. Note that for each $i$, $\dim(C_{i}) =
dim(T)$.

\bigbreak
\noindent
\emph{Claim 4.} For each $i$ there is a surjective
holomorphic map from $C_{i}$ to $T$ (so
finite-to-one outside a proper Zariski closed subset
$D_{i}$ of $C_{i}$).

\begin{proof} By Step I we have a surjective holomorphic map
$\pi:G^{*}\rightarrow {\mathbb P}^{1}$ such that $\pi^{-1}(H) =
G$ and
$\pi|G$ is precisely the canonical surjective homomorphism
from $G$ to $H$. So $G^{*}\setminus G$ is $\pi^{-1}(S)$
where $S = {\mathbb P}^{1}\setminus H$. Consider the map
$\mu(g,h) = g\cdot h^{-1}$ from $G\times_{H}G$ to $T$. This
is definable and holomorphic so extends to a meromorphic map
from $G^{*}\times_{{\mathbb P}^{1}}G^{*}$ to $T$, which we also
call $\mu$. By Lemma 3.3 of \cite{Fujiki} this map is
holomorphic. In particular, for any $C_{i}$ and $x\in C_{i}$,
$\mu(-,x)|C_{i}$ is a holomorphic map from $C_{i}$ into $T$.
We must show that for suitable $x\in C_{i}$, this is
surjective.

For $g\in G$, let $f_{g}$ be the
meromorphic map from
$G^{*}$ to $G^{*}$ whose restriction to $G$ is
multiplication by $g$ (so $f_{g}$ is ${\bar \tau_{g}}$ in
previous notation).

Let $t\in T$. The map taking $x\in G$ to
$\mu(t \cdot x,x)\in T$ is the constant map with value $t$. It
follows that whenever $x\in G^{*}$ and $f_{t}$ is single
valued at $x$ then $\mu(f_{t}(x),x) = t$. Choose $x_{0}$
generic in $C_{i}$. Then for a dense open set $V$ of $t's$
in $T$, $f_{t}(x_{0})$ is defined and in
$C_{i}$. So for each $t\in V$, $\mu(f_{t}(x_{0}),x_{0}) = t$.
Thus $\mu(-,x_{0})|C_{i}:C_{i}\rightarrow T$ is generically
surjective, so surjective. In any case this map is finite.
\end{proof}

\noindent
  By Lemma 3.3, for each
$i$, $f_{g}$ induces a generic holomorphic action of $G$ on
$C_{i}$. Let $K_{i}$ be the subgroup of $G$ consisting of
those
$g\in G$ such that for all $x$ in some Zariski open subset
of $C_{i}$, $f_{g}(x) = x$. For dimension reasons $K_{i}$ is
a definable infinite subgroup of $G$, so contains $T$.
Moreover we have an induced generic action of $G/K_{i}$ on
$C_{i}$. Let $D_{i}$ be as in Claim 4. Let
$x_{0}\in C_{i}$ be such that for any $h$ in some dense
Zariski open subset of $G/K_{i}$, $h.x_{0} \in
C_{i}\setminus D_{i}$. This gives a meromorphic map from
$G/K_{i}$ to $C_{i}$ whose image contains infinitely many
points outside $D_{i}$. Composing with the holomorphic map
from $C_{i}$ to $T$ given by Claim 4 yields a meromorphic
nonconstant map from ${\mathbb P}^{1}$ into $T$ which is
impossible.

Thus $K_{i} = G$ and the generic action of $G$ on $C_{i}$ is
trivial. This holds for each $i$. So the generic action of
$G$ on $G^{*}\setminus G$ is trivial. That is, if $g\in G$,
then the meromorphic mapping $f_{g}:G^{*} \rightarrow
(G^{*}\setminus G)$ agrees with the identity map on a dense
Zariski open set. This implies that $f_{g}$ is the identity
map. So $G^{*}$ witnesses $G$ being Fujiki meromorphic. The
proof of Lemma 3.8 is complete as well as the proofs of
Theorem 3.1 and Corollary 3.2
\end{proof}

\section{Composition series}
In this section we will prove
\begin{Theorem} Suppose $G$ is a connected meromorphic
group. Then $G$ is regular (in the sense of Definition 2.6).
\end{Theorem}

\bigbreak
We first state a
consequence of Theorem 3.1

\begin{Proposition} Suppose the connected meromorphic group
$G$ is simple, in the sense that $G$ has no
nontrivial, connected normal definable subgroup. Then $G$ is
either
\begin{itemize}
\item[(i)] a (almost simple) noncommutative algebraic group,
\item[(ii)] ${\mathbb G}_a$ or ${\mathbb G}_m$,
\item[(iii)] a simple abelian variety, or
\item[(iv)] a strongly minimal modular complex torus.
\end{itemize}
\end{Proposition}
\begin{proof} Simplicity of $G$ together with the
the dichotomy theorem from
\cite{Hrushovski-Zilber2}, implies that
$G$ is either nonorthogonal to
${\mathbb P}^{1}$ (namely has nontrivial algebraic reduction) or
$G$ is modular. In the first case, $G$ is an algebraic
group, so (i), (ii) or (iii) hold. In the second case, every
definable subset of $G$ is a Boolean combination of cosets
of definable subgroups. Simplicity implies that $G$ is
strongly minimal. By 3.2, $G$ is a complex torus.
\end{proof}

The following will be crucial. The special case for
Fujiki-meromorphic groups in the class ${\mathcal C}$ was
proved in \cite{Fujiki}. In any case the ``classical" theory
of groups of finite Morley rank enters the picture.
\begin{Lemma} Let $1 \rightarrow L \rightarrow G \rightarrow
H \rightarrow 1$ be an exact sequence of connected
meromorphic groups and suppose that $L$ and $H$ are
(meromorphically isomorphic to) linear algebraic groups.
Then so is $G$.
\end{Lemma}
\begin{proof}
Note that if $G$ satisfies the hypotheses of
the lemma and $G_{1}$ is a connected definable subgroup of
$G$, or an image of $G$ under a meromorphic homomorphism then
$G_{1}$ satisfies the hypotheses too (for suitable
$L_{1}, H_{1}$).

We will prove the lemma by
induction on $\dim(G) = n$. We consider various possibilities
for $G$,

\noindent
\emph{Case I.} $G$ has an an infinite center.

By the hypotheses, $Z(G)$ contains an infinite definable
linear algebraic group and thus, by the structure of
commutative linear algebraic groups, $Z(G)$ contains a
definable $1$-dimensional connected linear algebraic group
$A$. $A$ is normal in $G$, so by induction hypothesis $G/A$
is (meromorphically isomorphic to) a linear algebraic
group. It makes sense to talk about the algebraic dimension
$a(G)$ of $G$. Note that $\dim(G/A) = n-1$ so the map
$G\rightarrow G/A$ witnesses that $a(G) \geq n-1$. If $a(G)
= n$, then $G$ is already isomorphic to an algebraic group,
so a linear algebraic group. Otherwise $a(G) = n-1$, and it
is well-known (see \cite{Ueno}) that the general fiber of the
algebraic reduction $\pi:G\rightarrow X$ is an elliptic
curve $E$. But the map $G \rightarrow G/A$ must
meromorphically factor through
$\pi$. The general fiber of the first map is ${\mathbb P}^{1}$
and thus we see that $E$ is an image of ${\mathbb P}^{1}$ under
a meromorphic (i.e. rational) map, a contradiction. Thus $G$
is linear algebraic.

\bigbreak
\noindent
{\em Case II.} $G$ is solvable.

We may assume, by Case I, that $Z(G)$ is finite. But then
$G/Z(G)$ is centerless and easily $G$ is linear algebraic
iff $G/Z(G)$ is. So we may assume that $G$ is centerless.
As is shown in
Chapter 3 of \cite{Poizat} or Chapter 9 of
\cite{Borovik-Nesin}, the commutator subgroup
$G'$ of
$G$ is connected, and nilpotent, so $Z(G')$ is infinite and
contains a minimal definable connected $G$-normal subgroup
$A$. $G/G'$ defines an infinite group of automorphisms of
$A$. By again results in \cite{Poizat} or
\cite{Borovik-Nesin},
$A$ is the additive group of a definable field $K$. As $A$
is by hypothesis linear algebraic, the field $K$ has to be
(definably isomorphic to) ${\mathbb C}$ and
$\dim(A) = 1$. $G/A$ is by induction hypothesis algebraic,
and as in Case 1 we deduce that $G$ is (definably) algebraic.

\bigbreak
\noindent
{\em Case III.} $G$ is nonsolvable.
\newline
Note that $G$ is among other things a connected complex Lie
group, and as such we have the Levi-Malcev decomposition $G
= R.S$ where $R$ is the maximal normal solvable connected
subgroup of $G$, $S$ is a semisimple (complex) Lie
group (unique up to conjugacy in $G$), and $R\cap S$ is
discrete. Finite Morley rank considerations (see 5.38 in
\cite{Borovik-Nesin}) show that
$R$ is definable, so linear algebraic by Case II (or
induction). It is probably then well-known that $G$ must be
isomorphic as a complex Lie group to a linear algebraic
group. However we want $G$ to be definably isomorphic to
such a group. So we must do a little more work, although
maybe there is a more direct way. We will show

\noindent
{\em Claim. $S$ is a definable subgroup of $G$.}

\begin{proof} We may assume that $R$
is a proper, nontrivial subgroup of $G$, definably
isomorphic to a linear algebraic group. We will first reduce
to the case where $R$ is commutative and unipotent. Let
$H$ be the connected component of the center of the
commutator subgroup of $R$. $H$ is then a nontrivial
commutative connected linear algebraic group, normal in
$G$. So $H$ is the direct product $U.T$ of a commutative
unipotent group
$U$ and an algebraic torus $T$. Note that both $U$ and $T$
are definable connected normal subgroups of $G$. $T$ has no
infinite definable group of automorphisms, so is central in
$G$. By Case I we may assume $T$ to be trivial. Thus $H = U$
is unipotent. By induction hypothesis, $G/H$ is linear
algebraic. Clearly $R/H$ is the solvable radical of $G/H$.
Thus $G/H$ is an almost direct product of $R/H$ with a
semisimple algebraic group $G_{1}/H$ (where $G_{1}$ is a
definable connected subgroup of $G$ containing $H$). As $S$
is unique up to conjugacy we may assume that $G_{1} = H.S$.
Note that the homomorphism $\mu: G_{1}\rightarrow
G_{1}/H$ is an isomorphism on $S$.

Note that $G_{1}/H$ is linear algebraic, by induction
hypothesis among other things.  Now $S$ (being semisimple) is
isomorphic (uniquely) as a complex Lie group to a linear
algebraic group, so it makes sense to talk about an element
of $S$ being unipotent.
$S$ is an almost direct product of almost simple groups
$S_{1},.,S_{r}$. Fix a nontrivial unipotent element $a$ in
some $S_{i}\setminus H$. Work now inside the definable group
$G_{1}$. Let $a_{1} = \mu(a) \in G_{1}/H$. $a_{1}$ is then
unipotent, and let $U_{1}$ be a $1$-dimensional definable
unipotent subgroup of $G_{1}/H$ containing $a_{1}$. Let
$U_{2} =
\mu^{-1}(U_{1})$. Then $U_{2}$ is an extension of a
unipotent linear algebraic group ($U_{1}$) by a linear
algebraic unipotent group $H$ so is (by induction) linear
algebraic unipotent. $a\in U_{2}$. We can find a definable
commutative connected subgroup $U_{3}$ of $U_{2}$ containing
$a$. $U_{3}$ is definably a vector space over ${\bf C}$. The
$1$-dimensional subspace $U_{4}$ generated by $a$ is a
definable subgroup of $G$ contained in $S_{i}$. We have
found an infinite  connected subgroup $U_{4}$ of
$S_{i}$ which is definable in $G$. The group generated by
all the $U_{4}^{g}$, where $g\in S_{i}$ is definable and
must be equal to $S_{i}$. So $S_{i}$ is definable. As $i$
was arbitrary, $S$ is definable. The claim is proved. \end{proof}

\vspace{2mm}
\noindent
We want $S$ to be {\em definably} isomorphic to a linear
algebraic group. Recall that $S$ as a complex Lie group is
the almost direct product of almost simple (discrete centre)
groups $S_{1},..,S_{r}$. As the centre of $S$ is definable,
each $S_{i}$ has finite centre. By considering centralizers,
each $S_{i}$ is definable. $S_{i}$ being almost simple is
almost strongly minimal, hence by the validity of the Zilber
trichotomy in ${\cal A}$, modular or nonorthogonal to
${\mathbb P}^{1}$. But $S_{i}$ is nonabelian. So it must be
nonorthogonal to ${\mathbb P}^{1}$, so algebraic. Thus $S$
is definably an algebraic group, so by semismplicity, linear
algebraic.
\newline
Finally $G$, being the almost semidirect product of linear
algebraic $R$ with linear algebraic $S$, must be linear
algebraic. Case III is complete, as well as the proof of
Lemma 4.3.

\end{proof}

{\bf Proof of Theorem 4.1.}
The proof will by induction on $\dim(G)$.

We first deal with
the case when $G$ is commutative.
Let $H$ be a minimal definable connected subgroup of $G$.
By 4.2, $H$ is either (a) a linear algebraic group or (b) a
simple complex torus. If $G = H$ we are finished. Otherwise,
applying the induction hypotheses, $G/H$ is definably an
extension of a complex torus $T$ by a linear algebraic
group $L/H$. In case (a) by Lemma 4.3, $L$ is linear
algebraic. So
$G$ is definably an extension of $T$ by $L$, and we finish.
So suppose (b) holds. If $L/H$ is trivial, $G$ is an
extension of a complex torus by a complex torus, so also a
complex torus. Otherwise let $L_{1}/H$ be a $1$-dimensional
subgroup of
$L/H$. Then $L$ is definably an extension of ${\mathbb G}_a$ or
${\mathbb G}_m$ by the simple complex torus $T$. By Corollary 3.2
(ii)
$L$ splits, yielding a $1$-dimensional linear algebraic
subgroup of $G$. We are now back in case (a). This proves
Case I.

\bigbreak
We now deal with the general case.

If $G$ has no proper normal nontrivial definable connected
subgroup, then we are finished by 4.2

Otherwise let $H$ be a proper normal nontrivial connected
subgroup of $G$. The induction hypothesis applies to $H$. If
the maximal connected linear algebraic normal subgroup $L$
of $H$ is nontrivial, then $L$ is normal in $G$
as $L$ is characteristic in $H$ by Remark 2.7(i), and by
applying the induction hypothesis to $G/L$ and applying
Lemma 4.3 we finish. Otherwise $H$ is a complex torus. As a
complex torus has no infinite definable group of
automorphisms, $H$ is central in $G$. The induction
hypothesis applies to $G/H$. If the latter is a complex
torus so is $G$. Otherwise $G$ has a connected normal
definable subgroup $G_{1}$ containing $H$ such that
$G_{1}/H$ is linear algebraic. If $G_{1}/H$ is semisimple
(equivalently contains no infinite normal solvable
subgroup), then $G_{1}$ is the almost direct product of
its commutator subgroup $G_{1}'$ and $H$. $G_{1}'$ is
semisimple so (definably) linear algebraic. We conclude by
applying the induction hypothesis to $G/G_{1}'$ and using
Lemma 4.3.

If $G_{1}/H$ is not semisimple, then there is a definable
nontrivial connected subgroup $A$ of $G_{1}$ containing $H$
such that $A$ is normal in $G$ and $A/H$ is commutative (and
linear algebraic). $A/H$ is (definably) a product of
${\mathbb G}_a$'s and ${\mathbb G}_m$'s. Remember that $H$ is central in
$A$, so if $A$ is not commutative then the commutator map
yields a nonconstant meromorphic map from $A/H \times A/H$
into the complex torus $H$, which is impossible. So $A$ has
to be commutative. By induction hypothesis, or by the first
part of the proof (if $A = G$), $A$ has a definable
connected (nontrivial) subgroup which is normal in $G$ and
(definably) linear algebraic. By induction hypothesis and
Lemma 4.3 we finish.

This completes the proof of Theorem 4.1. \hspace{\fill} $\dashv$

\section{Additional remarks and questions}
The following was first proved by the second author
\cite{Scanlon}. We give a quick proof using Theorem 1.1.
\begin{Proposition} Let $X$ be a strongly minimal compact
complex manifold. Then $X$ is either a (smooth projective)
algebraic curve, a complex torus, or is trivial.
\end{Proposition}
\begin{proof}
  Suppose that $X$ is neither trivial, nor an
algebraic curve. Then $\dim(X) > 1$ and by
\cite{Hrushovski-Zilber1} and
\cite{Hrushovski-Zilber2}, there is a strongly minimal group
$G$ definable in $X$. By Theorem 1.1 $G$ must be a complex torus,
nonorthogonal to $X$. Nonorthogonality is witnessed by an
analytic subset
$\Gamma$ of $X\times A$ which projects generically
finite-to-one on each of $X$ and
$G$. As both
$X,G$ are strongly minimal, both projections $\pi_{1}:\Gamma
\rightarrow X$ and $\pi_{2}:\Gamma \rightarrow G$ are
finite-to-one, and $\Gamma$ is strongly minimal with
$\dim(\Gamma) = \dim(X) = \dim(G) > 1$. Replacing $\Gamma$
by its normalization we may assume
$\Gamma$ is normal. $\Gamma$ has no proper infinite analytic
subsets, in particular no co-dimension $1$ analytic subsets.
By the purity of branch theorem, $\pi_{2}$
is an unramified covering. Thus $\Gamma$ is a complex torus.
As $\pi_{1}$ is also an unramified covering $X$ is a complex
torus.
\end{proof}

We can more generally give  natural necessary and sufficient
conditions for a type of
$U$-rank $1$ to be trivial. First we call a compact complex manifold
(or space) $X$ simple if there is no definable family
$\{Y_{t}:t\in T\}$ of positive dimensional proper analytic
subvarieties of $X$ such that
$\bigcup\{Y_{t}:t\in T\}$ contains a Zariski open subset of
$X$. It is easy to see that $X$ is simple if and only if its
generic type $p_{X}$ has $U$-rank $1$. A Kummer manifold is a
compact complex space which is bimeromorphic with a space
of the form $T/G$ where $T$ is a complex torus and $G$ a
finite group of (holomorphic) automorphisms of $T$.

The following proposition explains the
trichotomy within simple compact complex spaces between algebraic
curves, nonalgebraic Kummer manifolds, and manifolds of zero
algebraic and Kummer dimension in terms of the Zilber trichotomy
for Zariski geometries.

\begin{Proposition} Suppose $X$ is a simple compact complex
manifold. Then
$p_{X}$ is trivial if and only if
\begin{itemize}
\item[(i)] the
algebraic dimension $a(X)$ of $X$ is $0$, and
\item[(ii)] there is
no surjective meromorphic map from $X$ to a Kummer manifold
(i.e. $k(X) = 0$ in the notation of \cite{Fujiki2}).
\end{itemize}
In particular, if $X\notin {\mathcal C}$ then $p_{X}$ is
trivial.
\end{Proposition}
\begin{proof}
($\Leftarrow$): suppose $p_{X}$ is
nonorthogonal to ${\mathbb P}^{1}$. Then clearly $a(X) >0$.
Suppose $p_{X}$ is nontrivial and modular. Then $p_{X}$ is
nonorthogonal to the generic type of a strongly minimal
modular torus $T$. Let
$a$ be a generic point of $X'$ (i.e. a realization of
$p_{X}$). Then there is generic $b\in T'$ in $acl(a)$. Let
$\{b_{1}, \ldots, b_{n}\}$ be the finite set of realizations of
$tp(b/a)$. Then by the modularity of $T$,
$(b_{1}, \ldots, b_{n})$ is a generic point of a translate of a
(strongly minimal) subtorus $S$ of $T^{n}$. After
translating we may assume that $(b_{1}, \ldots , b_{n})$ is a
generic point of $S$. By elimination of imaginaries the
finite set $\{b_{1}, \ldots, b_{n}\}$ is coded by some $c$. As
$c\in dcl(a)$, we obtain a surjective meromorphic map $f$
from $X$ to a compact complex manifold $Y$ where $c$ is a
generic point of $Y$. But the map taking $(b_{1}, \cdots, b_{n})$ to
$c$ extends to a meromorphic map from $S$ to $Y$. Modularity
of $S$ implies that this induces a bimeromorphic map between
$S/G$ and $Y$ for some finite group $G$ of automorphisms of
$S$. Thus $Y$ is Kummer. So if (i) and (ii) hold, the only
possibility left for $p_{X}$ is to be trivial.

\noindent
($\Rightarrow$):  If (i) fails then $p_{X}$ is nonorthogonal
to ${\mathbb P}^{1}$ so is nontrivial. If (ii) fails and there is
a surjective meromorphic map to the Kummer manifold $Y$ then
clearly $Y$ is also simple and its generic type is
nontrivial. So $p_{X}$ is nontrivial.

The ``in particular'' clause follows from the
observations that any Kummer manifold $Y$ is in ${\mathcal C}$ as well as
any
compact complex manifold which maps meromorphically and
generically finite-to-one on $Y$.
\end{proof}

Our classification of meromorphic groups together with
Fact 2.8 (i) shows that any meromorphic group $G$ is of
``type ${\mathcal C}$", that is already definable in the
structure ${\mathcal C}$. As ${\mathcal C}$ is saturated in a
suitable language,  our results also classify definable
groups in all models of
$Th({\mathcal C})$. What about groups definable in elementary
extensions of ${\mathcal A}$? Here is a conjecture:
\begin{Conjecture} Let $G$ be a definable group in ${\mathcal
A}'$. Then $G$ is definably (in ${\mathcal A}'$) isomorphic to a
group $H$ definable in the reduct ${\mathcal C}'$.
\end{Conjecture}

The conjecture can be restated in terms of families of
groups in ${\mathcal A}$.

\end{document}